\documentclass[11pt]{amsproc}

\usepackage{amsmath, amsthm, amssymb, amsfonts, amscd, graphicx, endnotes, tikz, color, mathrsfs, float, array}

\usepackage{mathtools}

\usepackage[margin=01.4in]{geometry}

\usetikzlibrary{arrows}

\usepackage{hyperref}

\usepackage[nodisplayskipstretch]{setspace}

\def\NN{{\mathbb N}}

\def\QQ{{\mathbb Q}}
\def\PP{{\mathbb P}}
\def\QQ{{\mathbb Q}}
\def\RR{{\mathbb R}}

\def\ZZ{{\mathbb Z}}

\def\0{{\mathbf 0}}
\def\1{{\mathbf 1}}
\def\a{{\mathbf a}}
\def\b{{\mathbf b}}

\def\k{{\mathbf k}}

\def\w{{\mathbf w}}

\def\Escr{{\mathscr E}}

\def\Acal{{\mathcal A}}
\def\Bcal{{\mathcal B}}

\def\Xcal{{\mathcal X}}

\def\R{\mathrm{R}}

\def\min{\mathrm{min}}

\theoremstyle{plain}

\newtheorem*{thm-nolabel}{Theorem}

\newtheorem{thm}{Theorem}

\newtheorem{conj}{Conjecture}
\newtheorem{cor}[thm]{Corollary}
\newtheorem{prop}[thm]{Proposition}
\newtheorem{lem}[thm]{Lemma}

\theoremstyle{definition}
\newtheorem*{dfn}{Definition}

\newtheorem{ex}{Example}

\title[Arithmetic normality and pseudorandomness \\ of $p$-adic scaling maps]{Arithmetic normality and pseudorandomness \\ of $p$-adic scaling maps}

\author{Sage Bothwell}

\address{Sage Bothwell \newline Department of Mathematics \newline Oregon State University \newline Corvallis, Oregon, U.S.A. \newline bothwels@oregonstate.edu}

\author{Clayton Petsche}

\address{Clayton Petsche \newline Department of Mathematics \newline Oregon State University \newline Corvallis, Oregon, U.S.A. \newline petschec@oregonstate.edu}

\thanks{{\em Date:} November 27, 2026}

\begin{document}

\begin{abstract}
We consider the class of $p$-adic scaling maps $f:\ZZ_p\to\ZZ_p$.  Given such a map $f$ and a seed $\alpha\in\ZZ_p$, we consider the $p$-ary itinerary sequence associated to the forward $f$-orbit of $\alpha$ through the partition of $\ZZ_p$ into its $p$ distinct congruence classes modulo $p$.  We conjecture that these itinerary sequences are normal in the arithmetic case of a $p$-adic scaling polynomial $f(x)\in\QQ[x]$ with rational coefficients and a non-preperiodic rational seed $\alpha\in\QQ\cap\ZZ_p$, and we investigate numerically the plausibility of this conjecture in $2$-adic examples.  We also prove a weak normality result for arbitrary $p$-scaling maps.  We consider the itineraries of $2$-adic scaling maps as pseudorandom number generators, and we use using empirical statistical tests to assess their viability.
\end{abstract}

\maketitle

\tableofcontents

\newpage


\section{Introduction}\label{IntroSect}

\subsection{A pseudorandom binary sequence}  Define a sequence $\alpha_0,\alpha_1,\alpha_2,\alpha_3,\dots$ of integers by selecting any initial value $\alpha_0\in\ZZ$, and continuing recursively by 
\begin{equation}\label{LogisticRecursionIntro}
\alpha_{n+1}=\frac{\alpha_n(\alpha_n-1)}{2}.
\end{equation}
For each $\alpha\in\ZZ$, define $\bar{\alpha}=0$ if $\alpha$ is even and $\bar{\alpha}=1$ if $x$ is odd.  With $\alpha_0=5$, the first eight terms of $\alpha_n$ and $\bar{\alpha}_n$ are as follows.
\begin{center}
\begin{tabular}{lll}
$n$ & $\alpha_n$ &  $\bar{\alpha}_n$ \\
\hline
$0$ & 5 & 1  \\
$1$ & 10 &  0 \\
$2$ & 45 &  1 \\
$3$ & 990 &  0 \\
$4$ & 489555 & 1  \\
$5$ & 119831804235 &  1 \\
$6$ & 7179830653047766065495 &  1 \\
$7$ & 25774984103222155465711757501414912431764765 & 1 \\
\end{tabular}
\end{center}
Continuing further, the infinite sequence of bits $(\bar{\alpha}_0\bar{\alpha}_1\bar{\alpha}_2\bar{\alpha}_3\dots)$ begins 
\begin{equation}\label{BitStringExample}
1010111101101011001110110001100111011011001110110011001...
\end{equation}

The bit sequence (\ref{BitStringExample}) appears plausibly random, as if it were obtained using a series of flips of a fair coin.  But as a single sequence defined using a deterministic recursion, this sequence is not truly random in any meaningful sense.  Deterministic algorithms to define sequences which appear as if they were generated at random, known as pseudorandom number generators (PRNGs), have a wide variety of uses, for example in machine learning, cryptography, and numerical integration, among other fields \cite{AntunesHillPRNGML, KatzLindell, MR836749}.  

It should be noted that the bit sequence $(\bar{\alpha}_n)$ is far less expensive to calculate than it may at first appear, given the rapid growth of $(\alpha_n)$.  For each fixed $N\geq1$, if we want to calculate the first $N$ terms $\bar{\alpha}_0\bar{\alpha}_1\bar{\alpha}_2\dots \bar{\alpha}_{N-1}$, it suffices to consider instead the recursion 
\begin{equation*}
\beta_0=\alpha_0\pmod{2^N} \hskip1cm \beta_{n+1}=\frac{\beta_n(\beta_n-1)}{2} \pmod{2^N},
\end{equation*}
at each iteration choosing the least nonnegative residue modulo $2^N$.  In this situation (see Proposition \ref{NBitCalculationProp}) it holds that $\bar{\alpha}_n=\bar{\beta}_n$ for all $n=0,1,2,\dots,N-1$.  This means that arithmetic of $N$-bit integers is sufficient to calculate the length $N$ bit string $\bar{\alpha}_0\bar{\alpha}_1\bar{\alpha}_2\dots \bar{\alpha}_{N-1}$, a significant savings over the highly impractical calculation of the $2^n$-bit integers $\alpha_n$ for $0\leq n\leq N-1$. 

\subsection{$p$-adic scaling maps}  The first to propose sequences such as (\ref{BitStringExample}) as PRNGs were Woodcock-Smart \cite{woodcocksmart:chaosprng}, who were interested in $p$-adic dynamical analogues of the logistic family of maps in real dynamics.  Let $p$ be an arbitrary prime, and let $\ZZ_p$ denote the compact ring of $p$-adic integers, equipped with its usual absolute value $|\cdot|_p$.  Woodcock-Smart considered the dynamical system
\begin{equation}\label{pAdicLogisticIntro}
f:\ZZ_p\to\ZZ_p \hskip1cm f(x) = \frac{x^p-x}{p},
\end{equation}
which they called the $p$-adic logistic map in an analogy with the well-known family of logistic maps in real dynamics. The main result of \cite{woodcocksmart:chaosprng} is a proof that (\ref{pAdicLogisticIntro}) is conjugate to a Bernoulli shift on $p$ symbols.  

Kingsbery-Levin-Preygel-Silva \cite{MR2439398,MR2773138} later identified the wider class of {\em locally scaling} $p$-adic dynamical systems, which generalize the behavior of (\ref{pAdicLogisticIntro}).  Here we restrict attention to $p$-adic dynamical systems satisfying the following property, which represent a subclass of those considered in \cite{MR2439398,MR2773138}.

\begin{dfn}
A {\em $p$-adic scaling map} is a function $f:\ZZ_p\to\ZZ_p$ satisfying the property that $|f(x)-f(y)|_p=p|x-y|_p$ for all $x,y\in\ZZ_p$ with $|x-y|_p\leq1/p$.
\end{dfn} 

Examples of $p$-adic scaling maps $f:\ZZ_p\to\ZZ_p$ defined by explicit polynomials are quite plentiful and easily found.  In Proposition \ref{pScalingPolynomials} we show that if $u\in\ZZ_p$ satisfies $|u|_p=1$, and if $f$ differs from $\frac{u}{p}(x^p-x)$ by any polynomial in $\ZZ_p[x]$ of any degree, then $f(\ZZ_p)\subseteq\ZZ_p$ and $f:\ZZ_p\to\ZZ_p$ is a $p$-adic scaling map.  It is also possible to use Mahler expansions to describe a wide class of $p$-adic scaling maps which are not necessarily polynomials, as done by Kingsbery-Levin-Preygel-Silva \cite{MR2773138}. 

We denote the standard set of {\em $p$-ary digits} by
\[
\Sigma_p=\{0,1,2,\dots,p-1\}.
\] 
For each $\alpha\in\ZZ_p$, define $\bar{\alpha}$ to be the unique element of $\Sigma_p$ satisfying $\bar\alpha\equiv \alpha\pmod{p}$.  Thus 
\[
\ZZ_p\to\Sigma_p \hskip1cm \alpha\mapsto\bar\alpha
\]
is a surjective map whose fibers form the partition 
\begin{equation}\label{ZpPartition}
\ZZ_p=D_{1/p}(0)\cup D_{1/p}(1)\cup D_{1/p}(2)\cup\dots \cup D_{1/p}(p-1)\
\end{equation}
of $\ZZ_p$ into its $p$ disjoint discs of radius $1/p$.

Let $\NN=\{0,1,2,3,\dots\}$ be the set of nonnegative integers, and define
\[
\Sigma_p^\NN=\{\k=(k_0k_1k_2\dots)\mid k_n\in\Sigma_p\}
\]
to be the compact space of all infinite $p$-ary sequences.  Given any map $f:\ZZ_p\to\ZZ_p$ and any seed $\alpha\in\ZZ_p$, iteration of $f$ defines a $p$-ary sequence $(\overline{f^n(\alpha)})\in\Sigma_p^\NN$, which we refer to as the {\em itinerary} of the forward $f$-orbit of $\alpha$ through the partition (\ref{ZpPartition}).  The following theorem is due to Woodcock-Smart \cite{woodcocksmart:chaosprng} in the special case (\ref{pAdicLogisticIntro}), and generalized by Kingsbery-Levin-Preygel-Silva \cite{MR2773138}.

\begin{thm}[Woodcock-Smart \cite{woodcocksmart:chaosprng}, Kingsbery-Levin-Preygel-Silva \cite{MR2773138}]\label{TopConjThmIntro}
Let $f:\ZZ_p\to\ZZ_p$ be a $p$-adic scaling map.  The map 
\[
\ZZ_p\to\Sigma_p^\NN \hskip1cm \alpha\mapsto(\overline{f^n(\alpha)})
\] 
is a topological and measure preserving conjugacy of $f:\ZZ_p\to\ZZ_p$ with the shift map $\sigma:\Sigma_p^\NN\to\Sigma_p^\NN$ defined by $\sigma(k_0k_1k_2k_3\dots)=(k_1k_2k_3k_4\dots)$.
\end{thm}

\subsection{True randomness of $p$-adic scaling maps}  Although Woodcock-Smart and Kingsbery-Levin-Preygel-Silva did not choose to formulate their results in explicitly probabilistic terms, if we view $\ZZ_p$ as a probability space with respect to its normalized Haar measure, we can immediately obtain the following closely related result.  

\begin{thm}{\rm (True randomness of $p$-adic scaling maps)}\label{RandVarCorIntro}
Let $f:\ZZ_p\to\ZZ_p$ be a $p$-adic scaling map.  For each integer $n\geq0$ define 
\[
X_{f,n}:\ZZ_p\to\RR \hskip1cm X_{f,n}(\alpha)=\overline{f^n(\alpha)}.
\]
As a sequence of random variables on $\ZZ_p$, the $\{X_{f,n}\}_{n=0}^{\infty}$ are independent and identically distributed, taking equiprobable values in the set $\Sigma_p=\{0,1,2,\dots,p-1\}$. 
\end{thm}

In this paper we view Theorem \ref{RandVarCorIntro} as a probabilistic heuristic for arithmetic questions, and in particular we are primarily interested in $p$-adic scaling maps $f:\ZZ_p\to\ZZ_p$ defined by polynomials $f(x)\in\QQ[x]$ with rational coefficients, and with integer or rational seeds $\alpha\in\QQ$.  This would likely also be the primary case of interest in practical applications of $p$-adic scaling maps as PRNGs.  We point out that the statement of Theorem \ref{RandVarCorIntro} itself is of little practical use to study this situation, as this purely measure-theoretic result is insensitive to the behavior of any one seed or even any countable collection of seeds.

Perhaps the most intuitive property held (almost surely) by a random $p$-ary sequence is normality.  A $p$-ary sequence $\k=(k_n)\in\Sigma_p^\NN$ is {\em normal} if every finite word in the alphabet $\Sigma_p$ occurs in the sequence $\k$ with the expected asymptotic frequency of $1/p^{r}$, where $r$ is the length of the word.  (See $\S$\ref{NormalSeedsSect} for the more precisely stated definition.)  

\begin{dfn}
Let $f:\ZZ_p\to\ZZ_p$ be a $p$-adic scaling map.  A point $\alpha\in\ZZ_p$ is a {\em normal seed} for $f$ if its itinerary $(\overline{f^n(\alpha)})\in\Sigma_p^\NN$ is a normal sequence.
\end{dfn}

The following fact is immediate from Theorem \ref{TopConjThmIntro}.

\begin{cor}\label{ScalingAENormalIntro}
Let $f:\ZZ_p\to\ZZ_p$ be a $p$-adic scaling map.  Then Haar almost every $\alpha\in\ZZ_p$ is a normal seed for $f$.
\end{cor}

\subsection{Inverting the conjugacy and construction of an explicit normal seed}  Given a $p$-adic scaling map $f:\ZZ_p\to\ZZ_p$, the topological conjugacy 
\[
\ZZ_p\to\Sigma_p^\NN \hskip1cm \alpha\mapsto(\overline{f^n(\alpha)})
\] 
of Theorem \ref{TopConjThmIntro} is shown to be bijective using an indirect counting argument on the finite quotients $\ZZ_p/p^N\ZZ_p$.  Alternatively, we show in $\S$\ref{NormalSeedSect} that this conjugacy can be explicitly inverted using a recursion.  More precisely, given an arbitrary $p$-ary sequence $\k=(k_n)\in\Sigma_p^\NN$, we show how to construct the unique seed $\alpha\in\ZZ_p$ with itinerary $(\overline{f^n(\alpha)})$ equal to $\k$.  In particular, this result allows for the construction of an explicit normal seed $\alpha$ for $f$, written with arbitrary precision. 

\subsection{An arithmetic normality conjecture}  Recall that $\alpha\in\ZZ_p$ is called {\em $f$-preperiodic} if the forward orbit $\{f^n(\alpha)\mid n\geq0\}$ of $\alpha$ contains only finitely many distinct points.  It follows from the theory of heights in arithmetic dynamics (\cite{MR2316407} Thm. 3.12) that if $f$ is defined by a polynomial with coefficients $\QQ$, then there are only finitely many rational $f$-preperiodic points $\alpha\in\QQ$.

\begin{conj}[Arithmetic normality of $p$-adic scaling maps]\label{ScalingArithNormalConjIntro}
Let $f:\ZZ_p\to\ZZ_p$ be a $p$-adic scaling map defined by a polynomial $f(x)\in\QQ[x]$, and let $\alpha\in\QQ\cap\ZZ_p$.  If $\alpha$ is not $f$-preperiodic, then $\alpha$ is a normal seed for $f$.
\end{conj}

One way to view Conjecture \ref{ScalingArithNormalConjIntro} is as an arithmetic-dynamical unlikely intersection principle, as it predicts that two Haar-measure null subsets of $\ZZ_p$, namely the non-$f$-normal seeds and the rational seeds, intersect only in the finite set of rational $f$-preperiodic points.  This conjecture is arithmetic in the sense that it {\em cannot} be expected to hold for arbitrary $p$-adic scaling maps with the hypothesis $f(x)\in\QQ[x]$ omitted.  For example, if $f:\ZZ_p\to\ZZ_p$ is an arbitrary $p$-adic scaling map and $\alpha\in\ZZ_p$ is any non-preperiodic, non-normal seed for $f$ (these exist by the conjugacy Theorem \ref{TopConjThmIntro}), then $g(x)=f(x+\alpha)-\alpha$ is also $p$-adic scaling and 0 is a non-preperiodic, non-normal seed for $g$.

In $\S$\ref{CompNormalSect} we calculate the frequency ratios of all words $\w$ of length $\leq 3$ occurring in finite itinerary sequences $(\overline{f^n(\alpha)})$ associated to two sample examples of $2$-adic scaling maps defined by polynomials $f(x)\in\QQ[x]$.  The results of these calculations do seem to behave as expected according to Conjecture \ref{ScalingArithNormalConjIntro}; see Tables \ref{WordFreqTable1} and \ref{WordFreqTable2}.

A proof of Conjecture \ref{ScalingArithNormalConjIntro} seems far out of reach.  We are unable to provide even a single example of a $p$-adic scaling map $f:\ZZ_p\to\ZZ_p$, defined by a polynomial $f(x)\in\QQ[x]$, and having a single rational seed $\alpha\in\QQ\cap\ZZ_p$ which can be proved to be $f$-normal.  The situation is reminiscent of a well-known theorem due to Koksma \cite{MR1556918} in metric number theory, stating that for Lebesgue almost all $\theta>1$, the sequence $\{\theta^n\}_{n=0}^{\infty}$ is uniformly distributed modulo 1.  Despite this theorem, not a single naturally occurring example of such a parameter $\theta$ is known.  For example, while it seems plausible that $\{(3/2)^n\}_{n=0}^{\infty}$ is uniformly distributed modulo 1, there is no known proof or disproof of this statement.

As described above, we can use Proposition \ref{PrescribedItineraryProp} to construct an explicit normal seed $\alpha$ for $f$, written with arbitrary precision.  But there is no reason to suspect that the $\alpha$ so constructed is rational.  

Regarding the normality of integer seeds, we can prove the following statement, which unfortunately is quite a bit weaker than Conjecture \ref{ScalingArithNormalConjIntro}.  Given a $p$-adic scaling map $f:\ZZ_p\to\ZZ_p$ and an integer $R\geq1$, we say that $\alpha\in\ZZ_p$ is an {\em $R$-weakly normal seed} for $f$ if every $p$-ary word of length at most $R$ occurs at least once in the sequence $(\overline{f^{n}(\alpha)})$.  We can prove that most positive integers are $R$-weakly normal seeds for $f$ in an asymptotic sense.  Given a subset $S\subseteq\NN$, recall that the {\em density} of $S$ is defined to be
\[
\delta(S)=\lim_{M\to+\infty}\frac{|S\cap\{0,1,2,3,\dots,M-1\}|}{M},
\]
if the limit exists.

\begin{thm}\label{ScalingWeaklyNormalIntro}
Let $f:\ZZ_p\to\ZZ_p$ be a $p$-adic scaling map and let $R\geq1$ be an integer.  Then the set of all $\alpha\in\NN$ which are $R$-weakly normal seeds for $f$ has density 1 in $\NN$.  
\end{thm}

\subsection{$2$-adic scaling maps as PRNGs}  We now return to the question of using the itinerary sequences $(\overline{f^n(\alpha)})$ of $2$-adic scaling maps as PRNGs.  Woodcock-Smart \cite{woodcocksmart:chaosprng} focused on the case of individual seeds $\alpha\in\ZZ_2$, and rather than considering the itinerary sequence $(\overline{f^n(\alpha)})$ itself, they considered instead the approximations $(\overline{f^n_e(\alpha)})$ indexed by an integer parameter $e\geq1$.  Here, $f_e$ denotes the permutation of $\{1,2,3,\dots,2^e\}$ defined by $f_e(x)=\frac{1}{2}(x^2-x)\pmod{2^e}$, using the least positive congruence class representative modulo $2^e$.  The itinerary sequence $(\overline{f^n_e(\alpha)})$ is less expensive to calculate than $(\overline{f^n(\alpha)})$ for large $n$, and the two itineraries coincide for indices $n\leq e$.  A main focus of \cite{woodcocksmart:chaosprng} is to assess of pseudorandomness potential of $(\overline{f^n_e(\alpha)})$ for large $n>e$.

Our approach differs from \cite{woodcocksmart:chaosprng} in two ways.  First, we consider only the true itinerary sequences $(\overline{f^n(\alpha)})$ themselves, not the approximations $(\overline{f^n_e(\alpha)})$.

Second, rather than considering the itinerary sequence $(\overline{f^n(\alpha)})$ for any one particular individual seed $\alpha\in\ZZ_p$, we adopt the view commonly taken in the cryptography literature (e.g.~\cite{RosulekTJOC} $\S$5.1) that a PRNG is a string-lengthening function of a low entropy seed space.  As true randomness is difficult and expensive to obtain (commonly relying on the collection of data from a physical source of thermal noise), the purpose of a PRNG is the computationally inexpensive production of long pseudorandom bit strings, without the need for any additional entropy beyond the choice of seed.  

Let $f:\ZZ_2\to\ZZ_2$ be a $2$-adic scaling map defined by a polynomial $f(x)\in\QQ[x]$, and let $M<N$ be integers.  The PRNG associated to $f,M,N$ is the function
\begin{equation*}
G_f:\{0,1\}^M\to\{0,1\}^N
\end{equation*}
defined as follows.  Given a bit string $a_0a_1a_2\dots a_{M-1}\in\{0,1\}^M$, assemble it into the $2$-adic expansion of a nonnegative integer in the range $0\leq \alpha\leq 2^M-1$ by 
\[
\alpha=a_0+a_12+a_22^2+a_32^3+\dots+a_{M-1}2^{M-1},
\]
and define
\begin{equation}\label{ItineraryPRNGDef}
G_f(a_0a_1a_2\dots a_{M-1}) = b_0b_1b_2\dots b_{N-1} \hskip1cm b_n=\overline{f^n(\alpha)}.
\end{equation}

To illustrate with the same example as in (\ref{LogisticRecursionIntro}) and (\ref{BitStringExample}), we again select $f(x)=\frac{1}{2}(x^2-x)$, seed length $M=4$, and seed $a_0a_1a_2a_3=1010\in\{0,1\}^4$, which is the coefficient vector of the $2$-adic representation of
\[
\alpha=5=1\cdot2^0+0\cdot2^1+1\cdot2^2+0\cdot2^3.
\]
Selecting pseudorandom string length $N=12$, we have $G_f:\{0,1\}^4\to \{0,1\}^{12}$ and in particular we obtain from the seed $1010$ the pseudorandom string
\[
G_f(1010)=b_0b_1b_2b_3b_4b_5b_6b_7b_8b_9b_{10}b_{11}=101011110110
\]
where $b_n=\overline{f^n(5)}$, which is the first 12 bits of the itinerary sequence (\ref{BitStringExample}).

We consider the domain and codomain of the PRNG $G_f:\{0,1\}^M\to\{0,1\}^N$ as probability spaces with the standard uniform measures $\PP_M$ and $\PP_N$, respectively.  The length $M$ bit strings in the domain $\{0,1\}^M$ of $G_f$ are viewed as {\em random seeds}.  Arbitrary length $N$ strings in the codomain $\{0,1\}^N$ of $G_f$ are {\em true random strings}, and length $N$ strings in the image $G_f(\{0,1\}^M)$ of $G_f$ are {\em pseudorandom strings}.  Intuitively, we view $G_f$ favorably as a PRNG if the statistical properties of pseudorandom strings in $\{0,1\}^N$ closely resemble those of true random strings in $\{0,1\}^N$.

To assess the statistical strength of $G_f:\{0,1\}^M\to\{0,1\}^N$ as a PRNG, we propose to compare the distribution of word frequencies among true random strings and pseudorandom strings of length $N$, as follows.  

Given a word $\w\in\{0,1\}^r$ of some length $r\geq1$, and an integer $q=1,2,3,\dots$, let us denote by $\mathscr{E}_{q,\w}$ the set of strings in $\{0,1\}^N$ such that the number of occurrences of $\w$ in the string is within $q$ standard deviations of its expected value.  We compare the probabilities 

\begin{center}
\begin{tabular}{rm{9cm}}
$\PP_N(\mathscr{E}_{q,\w})=$  &  The probability that a random string $\b\in\{0,1\}^N$ satisfies $\b\in\mathscr{E}_{q,\w}$ \\
$\PP_M(G_f^{-1}(\mathscr{E}_{q,\w}))=$  & The probability that a random seed $\a\in\{0,1\}^M$ satisfies $G_f(\a)\in\mathscr{E}_{q,\w}$
\end{tabular}
\end{center}

\noindent
If
\begin{equation}\label{PRApproxComparison}
\PP_M(G_f^{-1}(\mathscr{E}_{q,\w}))\approx \PP_N(\mathscr{E}_{q,\w})\hskip1cm q=1,2,3,\dots,
\end{equation}
this suggests that the distribution of the frequency of the word $\w$ fails to meaningfully distinguish between pseudorandom and true random strings in $\{0,1\}^N$.  We may interpret this as one item of empirical evidence that the PRNG $G_f$ successfully produces low entropy, pseudorandom strings in $\{0,1\}^N$ which closely resemble high entropy, true random strings in $\{0,1\}^N$.  

In $\S$\ref{PRNGSect} we carry out the comparison (\ref{PRApproxComparison}) in sample examples, taking $M=12$ and $N=240$ and thus representing PRNGs which increase the lengths of bit strings by a factor of 20.  Our data considers the cases $q=1,2,3$, all words $\w\in\{0,1,00,01,10,11\}$ of length at most $2$, for two particular $2$-adic scaling maps.  The data, presented in Tables \ref{DistWordFreqTable1} and \ref{DistWordFreqTable2}, do seem to show a close resemblance between pseudorandom and true random word frequency distributions, supporting the viability of $G_f:\{0,1\}^M\to\{0,1\}^N$ as statistically sound PRNGs.

\subsection{Further discussion}  The idea of chaotic dynamical systems as a PRNGs goes back at least to a brief 1947 preliminary report of Ulam and von Neumann \cite{UlamVonNeumann}, who suggested the chaotic real logistic map $f:[0,1]\to[0,1]$ defined by  $f(x)=4x(1-x)$ for this purpose.  The idea was approached from a modern perspective by Phatak-Rao \cite{PhatakRao} in 1993.  A flurry of interest in real chaotic dynamical PRNGs followed (e.g. \cite{Andrecut, PareekPatidarSud}), but mostly failed to show strong practical advantages over conventional PRNGs.  Activity in this area has fallen out of fashion, although there have been some occasional new ideas in this direction (e.g. \cite{MachicaoBruno}).  

Real chaotic dynamical PRNGs have some disadvantages such as non-uniformity, as well as truncation error which can multiply rapidly given the expanding nature of chaotic dynamical systems. In contrast, the PRNGs $G_f$ considered here are expected to be statistically uniform, and can be precisely calculated at length $N$ working entirely within $\ZZ/p^N\ZZ$, avoiding the need for floating point arithmetic.  

In assessing the statistical qualities of PRNGs, there already exist well known standard tests such as the NIST Suite \cite{NIST} and Dieharder \cite{dieharder}; these are ad hoc collections of empirical statistical tests taking an individual bit sequence as input.  Rather than use these off the shelf, we opt instead for the method of empirical statistical assessment described above, which is designed to emphasize the nature of a PRNG not as a procedure for producing a single bit sequence, but rather as a string-lengthening function on a low entropy seed space. 

In this paper we do not consider the cryptographic security of the PRNGs associated to $p$-adic scaling maps; a PRNG $G:\{0,1\}^M\to\{0,1\}^N$ is considered {\em secure} if there does not exist a polynomial time algorithm to distinguish pseudorandom strings from true random strings (\cite{RosulekTJOC} $\S$5.1).  No secure PRNG is known to exist (such a discovery would imply that $\mathrm{P}\neq \mathrm{NP}$), but it would still be interesting to investigate whether the security of the PRNGs $G_f:\{0,1\}^M\to\{0,1\}^N$ can be reduced to well established hardness assumptions in complexity theory.

Finally, the $p$-adic scaling maps are also interesting examples in non-Archimedean Fatou-Julia theory, as their Julia sets (in the Berkovich sense) are equal to $\ZZ_p$, where their canonical measures coincide with Haar measure.  (This is worked out for $f(x)=\frac{1}{p}(x^p-x)$ in \cite{MR3890051} Ex. 5.30, but the same calculations apply more generally to arbitrary $p$-adic scaling maps).  As a consequence of this particularly simple example of a $p$-adic canonical measure, Petsche-Stacy \cite{MR3958063} used backward orbits of the $p$-adic scaling map $f(x)=\frac{1}{p}(x^p-x)$ to find the smallest known limit point of the height of totally $p$-adic algebraic numbers.

The authors would like to thank Mike Rosulek for helpful discussions.

Data reported in Tables \ref{WordFreqTable1}, \ref{WordFreqTable2}, \ref{DistWordFreqTable1}, and \ref{DistWordFreqTable2} were calculated in Matlab and SageMath, with code available at \href{https://github.com/SevenLanguages/pAdicPRNG}{github.com/SevenLanguages/pAdicPRNG}.  The combinatorial calculations supporting the proof of Lemma \ref{WordFreqLem} were assisted by GPT-6 and verified by the authors.


\section{True randomness of \texorpdfstring{$p$}{p}-adic scaling maps}

Let $p$ be a prime.  We denote by $\QQ_p$ the locally compact field of $p$-adic numbers, equipped with its absolute value $|\cdot|_p$ normalized as usual so that $|p|_p=1/p$.  For each $c\in\QQ_p$ and each $r\in\ZZ$, the disc centered at $c$ and with radius $1/p^r$ in $\QQ_p$ is denoted by
\begin{equation*}
\begin{split}
D_{1/p^r}(c) & =\{a\in\QQ_p\mid |a-c|_p\leq 1/p^r\}.
\end{split}
\end{equation*}

Let $\ZZ_p=D_1(0)=\{\alpha\in\QQ_p\mid|\alpha|\leq1\}$ denote the compact subring of $p$-adic integers.  Denoting the standard set of {\em $p$-ary digits} by $\Sigma_p=\{0,1,2,\dots,p-1\}$, recall that each $\alpha\in\ZZ_p$ has a unique expansion 
\[
\alpha=a_0+a_1p+a_2p^2+a_3p^3+\dots \hskip1cm a_m\in\Sigma_p
\]
which is convergent in the $p$-adic topology.  Define the map 
\begin{equation}\label{ReductionModp}
\ZZ_p\to\Sigma_p \hskip1cm \alpha\mapsto \bar{\alpha}:=a_0
\end{equation}
sending each $\alpha\in\ZZ_p$ to the unique $a_0\in\Sigma_p$ satisfying $\alpha\equiv a_0\pmod{p}$ in $\ZZ_p$.  

Denote by $\lambda_p$ the Haar measure on the compact group $\ZZ_p$, normalized so that $\lambda_p(\ZZ_p)=1$.  It follows that for each integer $r\geq0$ and center $c\in\ZZ_p$, the disc $D_{1/p^r}(c)$ has $\lambda_p$-measure equal to its radius $1/p^r$.

Let $\NN=\{0,1,2,3,\dots\}$ denote the set of nonnegative integers, and define
\[
\Sigma_p^\NN=\{\k=(k_0k_1k_2\dots)\mid k_n\in\Sigma_p\}
\]
to be the space of all infinite $p$-ary sequences in the alphabet $\Sigma_p$, viewed as a compact metric space with metric $d_p(\k,\k')=p^{-\min\{n\mid k_n\neq k'_n\}}$.  Let $\mu_p$ denote the Borel probability measure on $\Sigma_p^\NN$ determined by the property that, for each integer $N\geq1$ and each choice of $N$ $p$-ary digits $\ell_0,\ell_1,\ell_2,\dots,\ell_{N-1}\in\Sigma_p$, the cylinder set
\begin{equation}\label{CylinderSetMeasure}
\begin{split}
\{\k\in\Sigma_p^\NN\mid k_0=\ell_0,k_0=\ell_0,\dots,k_{N-1}=\ell_{N-1}\}
\end{split}
\end{equation}
has $\mu_p$-measure equal to $1/p^N$.

Theorem \ref{TopConjThmIntro} is a generalization of a result of Woodcock-Smart \cite{woodcocksmart:chaosprng}.  Kingsbery-Levin-Preygel-Silva \cite{MR2773138} have proved a further generalization of Theorem \ref{TopConjThmIntro} to locally scaling maps defined on more general compact-open subsets of local fields, using the framework of Markov shifts.  In this section, for completeness and to provide a more elementary approach, we include a proof of Theorem \ref{TopConjThmIntro}.

\begin{lem}\label{LogisticPreImageBaseCase}
Let $f:\ZZ_p\to\ZZ_p$ be a $p$-adic scaling map.  Then for each $k\in\Sigma_p$, the map $f$ restricts to a bijection from the disc $D_{1/p}(k)$ onto $\ZZ_p$.
\end{lem}

\begin{proof}
Fix $k\in\Sigma_p$ and denote by $f_k:D_{1/p}(k)\to\ZZ_p$ the restriction of $f$ to $D_{1/p}(k)$.  Given distinct $a,b\in D_{1/p}(k)$, we have $|a-b|_p\leq1/p$ and hence using the scaling property of $f$ we have $|f(a)-f(b)|_p=p|a-b|_p\neq0$, and therefore $f_k$ is injective.  

To see that $f_k$ is surjective, let $N\geq1$ and for each $0\leq m\leq p^N-1$ define $a_m=k+pm\in D_{1/p}(k)$.  For distinct $0\leq m<n\leq p^N-1$ we have $|a_m-a_n|=|p(m-n)|_p\leq1/p$ and therefore using the scaling property of $f$ we have
\[
|f_k(a_m)-f_k(a_n)|_p=p|(a_m-a_n)|_p=p|p(m-n)|_p=|m-n|_p\geq1/p^{N-1}.
\]
Since the points $f_k(a_0),f_k(a_1),f_k(a_2),\dots,f_k(a_{p^N-1})$ are pairwise $p$-adically separated by at least $1/p^{N-1}$, they are distinct modulo $p^N$ and hence comprise a complete set of coset representatives for the quotient $\ZZ_p/p^N\ZZ_p$.  Since $N$ is arbitrary, this shows that $f_k(D_{1/p}(k))$ is dense in $\ZZ_p$.  Since $f_k$ is continuous and $D_{1/p}(k)$ is compact, $f_k(D_{1/p}(k))$ is closed, and hence $f_k(D_{1/p}(k))=\ZZ_p$, completing the proof.
\end{proof}

\begin{lem}\label{TrajectorySetIsADiscLemma}
Let $f:\ZZ_p\to\ZZ_p$ be a $p$-adic scaling map.  Let $N\geq0$, let $k_0,k_1,\dots,k_{N}\in \Sigma_p$, and define
\begin{equation*}
A_{k_0,k_1,k_2,\dots,k_{N}}=\{ a\in\ZZ_p\mid \overline{f^n(a)}=k_n\text{ for all } n=0,1,2,\dots,N\}.
\end{equation*}
Then $A_{k_0,k_1,k_2,\dots,k_{N}}$ is equal to a disc of radius $1/p^{N+1}$ in $\ZZ_p$.
\end{lem}

\begin{proof}
The $N=0$ case is trivial, as $A_{k_0}$ is just the fiber of $k_0$ under the map $\alpha\mapsto\bar{\alpha}$, which is equal to the disc $D_{1/p}(k_0)$.  For $N\geq1$, we will show that
\begin{equation}\label{TrajectorySetIsADiscPullbackIdentity}
A_{k_0,k_1,k_2,\dots,k_{N}}= D_{1/p}(k_0)\cap f^{-1}(A_{k_1,k_2,\dots, k_N}).
\end{equation}
Clearly $A_{k_0,k_1,k_2,\dots,k_{N}}\subseteq A_{k_0}=D_{1/p}(k_0)$, so to prove (\ref{TrajectorySetIsADiscPullbackIdentity}) it suffices to check that
\begin{equation}\label{TrajectorySetIsADiscPullbackIdentityDiscCase}
A_{k_0,k_1,k_2,\dots,k_{N}}= f^{-1}(A_{k_1,k_2,\dots, k_N}) \hskip1cm \text{whenever } a\in D_{1/p}(k_0).
\end{equation}
We have for $ a\in D_{1/p}(k_0)$ that
\begin{equation*}
\begin{split}
 a\in A_{k_0,k_1,k_2,\dots,k_{N}} & \Leftrightarrow \overline{f^{n}( a)}=k_n\text{ for all }n=1,2,\dots,N \\
	&  \Leftrightarrow f( a)\in A_{k_1,k_2,\dots,k_{N}} \\
	&  \Leftrightarrow  a\in f^{-1}(A_{k_1,k_2,\dots,k_{N}})
\end{split}
\end{equation*}
which proves (\ref{TrajectorySetIsADiscPullbackIdentityDiscCase}) and therefore (\ref{TrajectorySetIsADiscPullbackIdentity}).

We next prove that $A_{k_0,k_1,k_2,\dots,k_{N}}$ is nonempty via induction.  The $N=0$ case has already been established as $A_{k_0}=D_{1/p}(k_0)$.  If $N\geq1$ then $A_{k_0,k_1,k_2,\dots,k_{N}}$ is nonempty using the identity (\ref{TrajectorySetIsADiscPullbackIdentity}), because we may assume as an induction hypothesis that $A_{k_1,k_2,\dots,k_{N}}$ is nonempty, and the $f$-fiber of any point of $A_{k_1,k_2,\dots,k_{N}}$ meets the disc $D_{1/p}(k_0)$ by Lemma \ref{LogisticPreImageBaseCase}. 

Finally, given $ a\in A_{k_0,k_1,k_2,\dots,k_{N}}$ (which exists by the preceding paragraph), we will show that
\begin{equation}\label{TrajectorySetIsADiscGoal}
A_{k_0,k_1,k_2,\dots,k_{N}} = D_{1/p^{N+1}}( a)
\end{equation}
which will complete the proof of the lemma.  If $x\in D_{1/p^{N+1}}( a)$, thus $|x- a|_p\leq1/p^{N+1}$, then for all $0\leq n\leq N$ we may apply the scaling property of $f$ iteratively to obtain
\[
|f^n(x)-f^n( a)|_p=p^n|x- a|_p\leq p^n/p^{N+1}\leq 1/p
\]
and hence $\overline{f^n(x)}=\overline{f^n( a)}=k_n$ for all such $n$, verifying that $x\in A_{k_0,k_1,k_2,\dots,k_{N}}$.  

Conversely, if $x\in A_{k_0,k_1,k_2,\dots,k_{N}}$, then $\overline{f^n(x)}=k_n=\overline{f^n( a)}$ for all $n=0,1,2,\dots,N$ and therefore $|f^n(x)-f^n( a)|_p\leq 1/p$ for all such $n$.  We again apply the scaling property of $f$ to obtain $|f^N(x)-f^N( a)|_p=p^N|x- a|_p$ and therefore 
\[
\textstyle |x- a|_p=\frac{1}{p^N}|f^N(x)-f^N( a)|_p\leq\frac{1}{p^{N+1}}.
\] 
This shows that $x\in D_{1/p^{N+1}}( a)$, completing the proof of (\ref{TrajectorySetIsADiscGoal}) and the lemma.
\end{proof}

Lemma \ref{TrajectorySetIsADiscLemma} has the following simple consequence, which was mentioned in the introduction and which makes calculation of sequences such as (\ref{BitStringExample}) much faster than the calculation of actual forward orbits.  

\begin{prop}\label{NBitCalculationProp}
If $f:\ZZ_p\to\ZZ_p$ is a $p$-adic scaling map and if two seeds $\alpha,\beta\in\ZZ_p$ satisfy $\alpha\equiv \beta\pmod{p^{N+1}}$ for some $N\geq0$, then $\overline{f^n(\alpha)}=\overline{f^n(\beta)}$ for all $n=0,1,2,\dots,N$.  
\end{prop}
\begin{proof}
Fixing $\alpha\in\ZZ_p$, the set of all $x\in \ZZ_p$ satisfying $\overline{f^n(\alpha)}=\overline{f^n(x)}$ for all $n=0,1,2,\dots,N$ is a disc of radius $1/p^{N+1}$ containing $\alpha$, and since $\alpha\equiv \beta\pmod{p^{N+1}}$, $\beta$ is an element of this disc.
\end{proof}

\begin{proof}[Proof of Theorem \ref{TopConjThmIntro}]
Let $\pi:\ZZ_p\to\Sigma_p^\ZZ$ denote the map defined by $\pi(\alpha)=(\overline{f^n(\alpha)})$.  If $\k=(k_n)\in\Sigma_p^\NN$, then in the notation of Lemma \ref{TrajectorySetIsADiscLemma} we have that 
\[
\pi^{-1}(\k) = \bigcap_{N\geq0}A_{k_0,k_1,k_2,\dots,k_N}.
\]
This shows that $\pi^{-1}(\k)$ is a singleton set, since by Lemma \ref{TrajectorySetIsADiscLemma} this set is a nested intersection of discs of radius decreasing to zero in a complete metric space.  It follows that $\pi$ is a bijection.  For $\alpha\in\ZZ_p$ we have 
\[
\pi(f(\alpha))=(\overline{f^n(f(\alpha))})=(\overline{f^{n+1}(\alpha)})=\sigma((\overline{f^{n}(\alpha)}))=\sigma(\pi(\alpha))
\]
and therefore $\pi\circ f=\sigma\circ\pi$.

For each $N\geq0$, it follows from Lemma \ref{TrajectorySetIsADiscLemma} that $D\mapsto\pi^{-1}(D)$ is a bijection between the set of all $p^{N+1}$ discs of radius $1/p^{N+1}$ in $\Sigma_p^\NN$ onto the set of all $p^{N+1}$ discs of radius $1/p^{N+1}$ in $\ZZ_p$.  It follows that $\pi$ is an isometry, as in each of the spaces $\ZZ_p$ and $\Sigma_p^\NN$, the distance between two points is the radius of the smallest disc containing both points.  It also shows that $\pi_*(\lambda_p)=\mu_p$, as the measures $\lambda_p(\pi^{-1}(D))$ and $\mu_p(D)$ coincide for all discs $D$ in $\Sigma_p^\NN$, which form a base of open sets for the Borel $\sigma$-algebra.
\end{proof}

\begin{proof}[Proof of Theorem \ref{RandVarCorIntro}]
We first show that
\begin{equation}\label{IdentDistProb}
\lambda_p(\{ a\in\ZZ_p\mid X_{f,n}( a)=k\})=1/p.
\end{equation}
for each $n\geq0$ and $k=0,1,2,\dots,p-1$, which verifies that the $X_{f,n}$ are identically distributed and take equiprobable values in the set $\Sigma_p$.  In the notation of Lemma \ref{TrajectorySetIsADiscLemma} we have a partition
\[
\{ a\in\ZZ_p\mid X_{f,n}( a)=k\}=\bigcup_{k_0,k_1,\dots,k_{n-1}\in\Sigma_p}A_{k_0,k_1,\dots,k_{n-1},k}.
\]
By Lemma \ref{TrajectorySetIsADiscLemma}, this partition shows that $\{ a\in\ZZ_p\mid X_n( a)=k\}$ is a disjoint union of $p^n$ discs of radius (and hence Haar measure) equal to $1/p^{n+1}$, which proves (\ref{IdentDistProb}).

To show that the $X_{f,n}$ are independent, we must show that for each $N\geq0$ and choice of $k_0,k_1,k_2,\dots,k_{N}\in\Sigma_p$, it holds that
\begin{equation}\label{IndependentProb}
\lambda_p\left(\bigcap_{0\leq n\leq N}\{  a\in\ZZ_p \mid X_{f,n}( a)=k_n\}\right)=\prod_{0\leq n\leq N}\lambda_p(\{ a\in\ZZ_p \mid X_{f,n}( a)=k_n\}).
\end{equation}
The right hand side of (\ref{IndependentProb}) is equal to $1/p^{N+1}$, by (\ref{IdentDistProb}).  The left hand side of (\ref{IndependentProb}) is also equal to $1/p^{N+1}$ because
\begin{equation*}
\bigcap_{0\leq n\leq N}\{ a\in\ZZ_p \mid X_{f,n}( a)=k_n\}=A_{k_0,k_1,k_2,\dots,k_N}
\end{equation*}
which is known to be a disc of radius $1/p^{N+1}$ from Lemma \ref{TrajectorySetIsADiscLemma}.
\end{proof}

The following result shows that $p$-adic scaling maps defined by polynomials can be easily constructed.

\begin{prop}\label{pScalingPolynomials}
Given $u\in\ZZ_p$ with $|u|_p=1$ and an arbitrary polynomial $f_0(x)\in\ZZ_p[x]$, define 
\begin{equation}\label{pScalingDegpForm}
f(x) = \frac{u}{p}(x^p-x) + f_0(x).
\end{equation}
Then $f(\ZZ_p)\subseteq\ZZ_p$ and $f:\ZZ_p\to\ZZ_p$ is a $p$-adic scaling map.  
\end{prop}

\begin{proof}
For $f(x)$ as in (\ref{pScalingDegpForm}), Fermat's theorem implies that $a^p-a\equiv 0\pmod{p}$ for all $a\in\ZZ_p$, from which it follows that $|\frac{u}{p}(a^p-a)|_p\leq1$.  Since $f_0(x)\in\ZZ_p[x]$, it follows that $f(a)\in\ZZ_p$ for all $a\in\ZZ_p$.  To see that $f:\ZZ_p\to\ZZ_p$ is a $p$-adic scaling map, for two variables $x$ and $y$ we may factor
\begin{equation*}
\begin{split}
f_0(x)-f_0(y) & =(x-y)G(x,y) \\
(x^p-x)-(y^p-y) & =(x-y)H(x,y)
\end{split}
\end{equation*}
for $G(x,y),H(x,y)\in\ZZ_p[x,y]$, and explicitly $H(x,y)=\sum_{i=0}^{p-1}x^iy^{p-1-i}-1$.  Therefore
\begin{equation}\label{DerivativeCalc}
\begin{split}
p(f(x)-f(y)) & = pf_0(x)-pf_0(y)+u(x^p-x)-u(y^p-y)  \\
	& =  (x-y)\left(pG(x,y)+uH(x,y)\right).
\end{split}
\end{equation}
Since $|p|_p=1/p$, the proof of the desired identity 
\begin{equation}\label{ExpandingGoal}
|f(a)-f(b)|_p = p|a-b|_p
\end{equation}
whenever $|a-b|_p\leq 1/p$ is complete if we can show that
\begin{equation}\label{ExpandingCalcConstant}
\begin{split}
\left|pG(a,b)+uH(a,b)\right|_p=1
\end{split}
\end{equation}
for all such $a,b\in\ZZ_p$.  Since $|a-b|_p\leq 1/p$ we have $a\equiv b\pmod{p}$ and therefore
\begin{equation}\label{ExpandingCalcCongruence}
\begin{split}
pG(a,b)+uH(a,b) & \equiv u\sum_{i=0}^{p-1}a^ib^{p-1-i}-u  \\
	& \equiv upa^{p-1}-u \\
	& \equiv-u \not\equiv0\pmod{p}.
\end{split}
\end{equation}
The congruence (\ref{ExpandingCalcCongruence}) implies (\ref{ExpandingCalcConstant}) and hence (\ref{ExpandingGoal}), completing the proof that $f(x)$ is a $p$-adic scaling map.
\end{proof}

\begin{ex}
The map
\begin{equation}\label{2adicScalingAlternateExample}
f:\ZZ_2\to\ZZ_2\hskip1cm f(x) = \frac{3}{2}x^2-\frac{1}{2}x+1.
\end{equation} 
is a $2$-adic scaling map according to Proposition \ref{pScalingPolynomials}, since $f(x)=\frac{1}{2}(x^2-x)+x^2+1$.  Since in practice one would normally choose integer seeds, an advantage of this map is that it has no integer preperiodic points.  This contrasts with $f(x)=\frac{1}{2}(x^2-x)$, for which $-2,-1,0,1,2,3$ are all $f$-preperiodic, and so these seeds generate eventually periodic itineraries which cannot exhibit pseudorandomness. 
\end{ex}


\section{A weak normality result for integer seeds}\label{NormalSeedsSect}

In this section we prove Theorem \ref{ScalingWeaklyNormalIntro} on weak normality for integer seeds.  We first require several definitions.  Let $\w=w_0w_1w_2\dots w_{r-1}\in\Sigma_p^r$ be a word of length $r\geq1$ in the alphabet $\Sigma_p=\{0,1,2,\dots,p-1\}$, and let $\k=(k_n)\in\Sigma_p^\NN$ be a $p$-ary sequence.  We say that {\em $\w$ occurs in $\k$} if there exists some index $n\geq0$ such that $k_{n+i}=w_i$ for all $0\leq i\leq r-1$; if needed, we may specify in this case that $\w$ {\em occurs in $\k$ at position $n$}.  

\begin{dfn}
A $p$-ary sequence $\k\in\Sigma_p^\NN$ is {\em normal} if, for each integer $r\geq1$ and each finite word $w\in\Sigma_p^r$ of length $r$, it holds that 
\[
\lim_{N\to+\infty}\frac{|\{0\leq n\leq N-1\mid \text{$\w$ occurs in $\k$ at position $n$}\}|}{N}=\frac{1}{p^r}.
\]
\end{dfn}

\begin{lem}
Let $\w\in\Sigma_p^r$ be a word of length $r\geq1$.  For each $N\geq0$, let $\Acal_N(\w)$ be the set of all $\k\in\Sigma_p^\NN$ such that $\w$ does not occur in $\k$ at position $n$ for all $n=0,1,2,\dots,N-1$.  Then $\lim_{N\to+\infty}\mu_p(\Acal_N(\w))=0$.
\end{lem}

\begin{proof}
Let $L=\lfloor N/r\rfloor$.  Then $(L-1)r\leq N-1$ and therefore $\Acal_N(\w)\subseteq \Bcal_N(\w)$, where $\Bcal_N(w)$ is the set of all $\k\in\Sigma_p^\NN$ such that $\w$ does not occur in $\k$ at position $n$ for all $n=0,r,2r,\dots,(L-1)r$.  Since the length-$r$ subwords of any sequence $\k$ beginning at positions $n=0,r,2r,\dots,(L-1)r$ are non-overlapping, the $\mu_p$-measure of $\Bcal_N(\w)$ is equal to $(1-p^{-r})^L$ using the cylinder set formula (\ref{CylinderSetMeasure}).  Therefore
\[
\mu_p(\Acal_N(\w))\leq\mu_p(\Bcal_N(\w)) =(1-p^{-r})^{\lfloor N/r\rfloor}\to0 \text{ as } N\to+\infty,
\]
completing the proof.
\end{proof}

\begin{prop}\label{SymbolicEquidistributionProp}
Let $\{\k_m\}_{m=0}^\infty$ be a $\mu_p$-equidistributed sequence of $p$-ary sequences $\k_m\in\Sigma_p^\NN$, and let $\w\in\Sigma_p^r$ be finite word of length $r$.  Then the set of all indices $m\in\NN$ such that $\w$ does not occur in $\k_m$ has density zero in $\NN$.
\end{prop}

\begin{proof}
Let $\Acal(\w)$ be the set of all $\k\in \Sigma_p^\NN$ such that $\w$ does not occur in $\k$; we must show that 
\[
\delta(\{m\in\NN\mid \k_m\in \Acal(\w)\})=0.
\]
For each $N\geq1$, let $\Acal_N(\w)$ be the set of all $\k\in\Sigma_p^\NN$ such that $\w$ does not occur in $\k$ at position $n$ for all $n=0,1,2,\dots,N-1$.  Thus $\Acal(\w)\subseteq \Acal_N(\w)$.  The set $\Acal_N(\w)$ is a finite union of cylinder sets of the form (\ref{CylinderSetMeasure}).  Since such cylinder sets are both open and closed in the space $\Sigma_p^\NN$, the characteristic function $\Xcal_{\Acal_N(\w)}:\Sigma_p^\NN\to\RR$ is continuous.  Using the $\mu_p$-equidistribution hypothesis of the sequence $\{\k_m\}_{m=0}^\infty$ we have 
\begin{equation}\label{Squeezed}
\begin{split}
\delta(\{m\in\NN\mid \k_m\in \Acal(\w)\}) & \leq \delta(\{m\in\NN\mid \k_m\in \Acal_N(\w)\}) \\
	& = \frac{1}{M}\sum_{m=0}^{M-1}\Xcal_{\Acal_N(\w)}(\k_m) \\
	& \to \int\Xcal_{\Acal_N(\w)}d\mu_p \\
	& = \mu_p(\Acal_N(\w))\to0 \text{ as }N\to+\infty,
\end{split}
\end{equation}
Since the density $\delta(\{m\in\NN\mid \k_m\in \Acal(\w)\})$ does not depend on $N$, it must be equal to zero by (\ref{Squeezed}), completing the proof.
\end{proof}

\begin{proof}[Proof of Theorem \ref{ScalingWeaklyNormalIntro}]
Recall that $f:\ZZ_p\to\ZZ_p$ is a $p$-adic scaling map.  By Theorem \ref{TopConjThmIntro}, the map $\Phi:\ZZ_p\to\Sigma_p^\NN$ defined by $\Phi(\alpha)=(f^n(\alpha))$ is a measure preserving isometry.  Since the sequence $\{0,1,2,3,\dots\}$ of nonnegative integers is Haar-equidsitributed in $\ZZ_p$, it follows that the sequence $\{\Phi(m)\}_{m=0}^\infty$ is $\mu_p$-equidistributed in $\Sigma_p^\NN$.  

Applying Proposition \ref{SymbolicEquidistributionProp} to the sequence $\{\k_m\}_{m=0}^{\infty}$ defined by $\k_m=\Phi(m)=(f^n(m))$, we deduce that if $\w$ is an arbitrary $p$-ary word, then 
\[
\delta(\{m\in\NN\mid \w \text{ does not occur in }(f^n(m))\})=0.
\]
For any $R\geq1$, define $W_R$ to be the set of all $p$-ary words of length at most $R$.  Since the union of any finite set of density zero subsets of $\NN$ still has density zero, we deduce
\[
\delta(\{m\in\NN\mid \text{There exists $\w\in W_R$ which does not occur in }(f^n(m))\})=0
\]
and therefore also its complementary statement
\[
\delta(\{m\in\NN\mid \text{Every $\w\in W_R$ occurs in }(f^n(m))\})=1,
\]
which is the desired conclusion of the Theorem.
\end{proof}


\section{Inverting the conjugacy and construction of a normal seed}\label{NormalSeedSect}

Given a $p$-adic scaling map $f:\ZZ_p\to\ZZ_p$ and an arbitrary $p$-ary sequence $\k=(k_n)\in\Sigma_p^\NN$, it follows from Theorem \ref{TopConjThmIntro} that there exists a unique $p$-adic seed $\alpha\in\ZZ_p$ with itinerary $(\overline{f^n}(\alpha))$ equal to $\k$.  In this section we show how to calculate $\alpha$ via an explicit recursion.  In particular, this result allows for the construction of an explicit normal seed $\alpha$ for $f$, written with arbitrary precision.  For simplicity we restrict attention to the $p$-adic scaling map (\ref{pAdicLogisticIntro}), but the methods readily generalize.

\begin{prop}\label{PrescribedItineraryProp}
Consider the $p$-adic scaling map $f:\ZZ_p\to\ZZ_p$ defined by $f(x)=\frac{1}{p}(x^p-x)$.  Let $\k=(k_n)\in\Sigma_p^\NN$ be any $p$-ary sequence, and suppose that
\[
\alpha=a_0+a_1p+a_2p^2+\dots\in\ZZ_p \hskip1cm a_m\in\Sigma_p
\]
is the unique seed for which $\overline{f^n(\alpha)}=k_n$ for all $n\geq0$.  For each $M\geq0$ write 
\begin{equation}\label{AlphaApproximated}
\alpha_M=a_0+a_1p+a_2p^2+\dots+a_Mp^M.
\end{equation}
Then
\begin{equation}\label{ItineraryRecurrence}
\begin{split}
a_0 & =k_0 \hskip1cm a_{m} = \overline{(-1)^m(k_m-f^m(\alpha_{m-1}))} \text{ for } m\geq1.
\end{split}
\end{equation}
\end{prop}

Since $\k=(k_n)$ is given and the $p$-adic number $\alpha_{m-1}$ depends only on the $p$-ary cofficients $a_0,a_1,a_2,\dots,a_{m-1}$, we may view (\ref{ItineraryRecurrence}) as a recurrence allowing for the calculation of the seed $\alpha$ to arbitrary precision.

As a sample application of Proposition \ref{PrescribedItineraryProp}, consider $f:\ZZ_2\to\ZZ_2$ defined by $f(x)=\frac{1}{2}(x^2-x)$ and consider the binary Champernowne sequence
\[
\k=(k_n)=(0\,1\,00\,01\,10\,11\,000\,001\,010\,011\,100\,101\,110\,111\,\dots)\in\Sigma_2^\NN,
\]
formed by listing all binary words of length $r=1$, followed by all all binary words of lengths $r=2$, and so on for $r=3,4,5,\dots$, concatenating each of these finite lists together into a single infinite binary string.  Since $\k$ is a normal sequence in $\Sigma_2^\NN$, the unique $\alpha\in\ZZ_2$ such that $(\overline{f^n(\alpha)})=\k$ is a normal seed for $f$.  Explicitly, to a $2$-adic precision of $25$ bits we have 
\[
\alpha=\sum_{n=0}^{\infty}a_n2^n\hskip1cm (a_n)=(0100010100111101100001001\dots).
\]

\begin{proof}[Proof of Proposition \ref{PrescribedItineraryProp}]
We will first show that if $1\leq n\leq N$ and $\alpha,h\in\ZZ_p$ with $|h|_p=1/p^N$, then 
\begin{equation}\label{DifferenceQuotientIdentity}
f^n(\alpha+h)=f^n(\alpha)+\frac{(-1)^nh}{p^{n-1}}\left(\frac{1}{p}+\theta\right)
\end{equation}
for some $\theta\in\ZZ_p$.  In the following calculation, each occurrence of $\theta_i$ for $i=1,2,3,,\dots$ represents an element of $\ZZ_p$, possibly depending on other parameters.

From (\ref{DerivativeCalc}), because $h\equiv0\pmod{p}$ we have
\begin{equation*}
\begin{split}
\frac{p}{h}(f(\alpha+h)-f(\alpha)) & =\sum_{m=0}^{p-1}(\alpha+h)^m\alpha^{p-1-m}-1 \\
	& \equiv p\alpha^{p-1}-1 \\
	& \equiv-1\pmod{p} 
\end{split}
\end{equation*}
which implies that
\begin{equation}\label{DifferenceQuotientIdentity1}
f(\alpha+h)=f(\alpha)-h\left(\frac{1}{p}+\theta_1\right);
\end{equation}
this is the $n=1$ case of (\ref{DifferenceQuotientIdentity}).  Assuming that $N\geq2$, we set $h'=-h(\frac{1}{p}+\theta_1)$ and note that $|h'|_p=|\frac{h}{p}|_p=1/p^{N-1}$.  We may therefore apply (\ref{DifferenceQuotientIdentity1}) with $f(\alpha)$ in place of $\alpha$ and $h'$ in place of $h$ to obtain
\begin{equation*}
\begin{split}
f^2(\alpha+h) & =f(f(\alpha)+h') \\
	& =f^2(\alpha)-h'\left(\frac{1}{p}+\theta_2\right) \\
	& =f^2(\alpha)+\frac{h}{p}\left(\frac{1}{p}+\theta_3\right) \\
\end{split}
\end{equation*}
which is the $n=2$ case of (\ref{DifferenceQuotientIdentity}).  More generally, for arbitrary $N\geq1$, we may continue this calculation iteratively for all $n=1,2,3,\dots,N$, resulting in the proof of (\ref{DifferenceQuotientIdentity}) for the full range $1\leq n\leq N$.

We turn to the proof of (\ref{ItineraryRecurrence}).  Trivially $a_0=\overline{\alpha}=\overline{f^0(\alpha)}=k_0$.  Recall that $\alpha=\sum_{m=0}^{\infty}a_mp^m\in\ZZ_p$ and $(k_n)\in\Sigma_p^\NN$ satisfy $(\overline{f^n(\alpha)})=(k_n)$.  With $\alpha_M$ defined as in (\ref{AlphaApproximated}), we have $|\alpha-\alpha_M|_p\leq1/p^{M+1}$, which by Lemma \ref{TrajectorySetIsADiscLemma} implies that the itineraries of $\alpha$ and $\alpha_M$ have the same first $M+1$ terms, that is $\overline{f^n(\alpha_M)}=\overline{f^n(\alpha)}$ for all $0\leq n\leq M$.  Taking $M=n$ we obtain
\[
\overline{f^n(\alpha_n)}=\overline{f^n(\alpha)}=k_n
\]
and therefore $ f^n(\alpha_n)\equiv k_n\pmod{p}$ for all $n\geq0$.  On the other hand, for $n\geq1$ we use (\ref{DifferenceQuotientIdentity}) with $h=a_np^n$ to obtain
\begin{equation*}
\begin{split}
f^n(\alpha_n) & =f^n(\alpha_{n-1}+a_np^n) \\
	& =f^n(\alpha_{n-1})+(-1)^na_n(1+p\theta)
\end{split}
\end{equation*}
and therefore 
\[
k_n \equiv f^n(\alpha_n) \equiv f^n(\alpha_{n-1})+(-1)^na_n\pmod{p}.
\]
We deduce that
\[
a_n\equiv(-1)^n(k_n-f^n(\alpha_{n-1}))\pmod{p},
\]
from which the desired identity (\ref{ItineraryRecurrence}) follows.
\end{proof}


\section{Assessing normality in numerical examples}\label{CompNormalSect}

In this section we assess the plausibility of Conjecture \ref{ScalingArithNormalConjIntro} by calculating several sample examples of itinerary sequences.  Let $f:\ZZ_2\to\ZZ_2$ be a $2$-adic scaling map, and let $\alpha\in\ZZ_2$ be a seed.  Given a binary word $\w=w_0w_1w_2\dots w_{r-1}\in\{0,1\}^r$ of length $r$ and a parameter $N\geq1$, define the $N$-th {\em word frequency ratio} of the itinerary $(\overline{f^n(\alpha)})$ by
\[
\R_{f,\alpha,N}(\w):=\frac{|\{0\leq m\leq N-1\mid \w\text{ occurs in $(\overline{f^n(\alpha)})$ at position $m$}\}|}{N}.
\]

According to the prediction of normality stated in Conjecture \ref{ScalingArithNormalConjIntro}, we should expect $\R_{f,\alpha,N}(\w)\approx1/2^r$ for large $N$.  

In Tables \ref{WordFreqTable1} and \ref{WordFreqTable2}, we calculate the word frequency ratios $\R_{f,\alpha,N}(\w)$ for itineraries of the two maps $f(x)=\frac{1}{2}(x^2-x)$ and $f(x)=\frac{3}{2}x^2-\frac{1}{2}x+1$, for $N=10000$, seeds $\alpha=5,6,7$, and all words of length $r=1,2,3$.  We see results consistent with Conjecture \ref{ScalingArithNormalConjIntro}, with deviations $|\R_{f,\alpha,N}(\w)-1/2^r|$ comparable to what one would expect for random bit sequences.  For example, examining the frequency ratio of the single bit $\w=1$ among the first $N=10000$ bits of a random bit sequence in $\{0,1\}^\NN$, a calculation shows a standard deviation of $\sigma=1/2\sqrt{N}=0.005$.  By the central limit theorem, about $68\%$ of randomly chosen long bit sequences have the property that their proportion of $1$s is within one standard deviation of the mean.  Examining Tables \ref{WordFreqTable1} and \ref{WordFreqTable2}, we see that 5 out of the 6 values of $\R_{f,\alpha,N}(\w)$ for $\w=1$ are within one standard deviation of the mean $1/2=0.5$.  We will give a more systematic formulation of this analysis in $\S$\ref{PRNGSect} concerning the pseudorandomness of itinerary sequences $(\overline{f^{n}(\alpha)})$.

\begin{table}[H]
\centering 
\caption{Word frequency ratios $\R_{f,\alpha,N}(\w)$ for itineraries of the map $f(x)=\frac{1}{2}(x^2-x)$, for $N=10000$, seeds $\alpha=5,6,7$, and all words of length $r=1,2,3$.}\label{WordFreqTable1}

\begin{tabular}{c|c|c|c|c}
$\w$  & $1/2^r$ & $\R_{f,5,N}(\w)$ & $\R_{f,6,N}(\w)$ & $\R_{f,7,N}(\w)$ \\
\hline
0     	& 0.5000 & 0.4982 & 0.4989 & 0.5004   \\
1     	& 0.5000 & 0.5018 & 0.5011 & 0.4996   \\
00   	& 0.2500 & 0.2474 & 0.2469 & 0.2483   \\
01   	& 0.2500 & 0.2508 & 0.2520 & 0.2521   \\
10  	& 0.2500 & 0.2509 & 0.2519 & 0.2522   \\
11    	& 0.2500 & 0.2509 & 0.2492 & 0.2474   \\
000  	& 0.1250 & 0.1216 & 0.1232 & 0.1218   \\
001  	& 0.1250 & 0.1258 & 0.1237 & 0.1265   \\
010  	& 0.1250 & 0.1217 & 0.1262 & 0.1264   \\
011  	& 0.1250 & 0.1291 & 0.1258 & 0.1257   \\
100  	& 0.1250 & 0.1258 & 0.1237 & 0.1265   \\
101  	& 0.1250 & 0.1251 & 0.1282 & 0.1257   \\
110  	& 0.1250 & 0.1291 & 0.1258 & 0.1258   \\
111  	& 0.1250 & 0.1218 & 0.1234 & 0.1216   \\
\end{tabular}
\end{table}

\begin{table}[H]
\centering 
\caption{Word frequency ratios $\R_{f,\alpha,N}(\w)$ for itineraries of the map $f(x)=\frac{3}{2}x^2-\frac{1}{2}x+1$, for $N=10000$, seeds $\alpha=5,6,7$, and all words of length $r=1,2,3$.}\label{WordFreqTable2}

\begin{tabular}{c|c|c|c|c}
$\w$  & $1/2^r$ & $\R_{f,5,N}(\w)$ & $\R_{f,6,N}(\w)$ & $\R_{f,7,N}(\w)$ \\
\hline
0     	& 0.5000 & 0.5033 & 0.5031 & 0.4927   \\
1     	& 0.5000 & 0.4967 & 0.4969 & 0.5073   \\
00   	& 0.2500 & 0.2524 & 0.2579 & 0.2416   \\
01   	& 0.2500 & 0.2509 & 0.2452 & 0.2511   \\
10  	& 0.2500 & 0.2510 & 0.2452 & 0.2511   \\
11    	& 0.2500 & 0.2457 & 0.2517 & 0.2562   \\
000  	& 0.1250 & 0.1240 & 0.1339 & 0.1194   \\
001  	& 0.1250 & 0.1284 & 0.1240 & 0.1222   \\
010  	& 0.1250 & 0.1274 & 0.1225 & 0.1212   \\
011  	& 0.1250 & 0.1235 & 0.1227 & 0.1299   \\
100  	& 0.1250 & 0.1284 & 0.1240 & 0.1222   \\
101  	& 0.1250 & 0.1226 & 0.1212 & 0.1289   \\
110  	& 0.1250 & 0.1235 & 0.1227 & 0.1300   \\
111  	& 0.1250 & 0.1222 & 0.1290 & 0.1262   \\
\end{tabular}
\end{table}


\section{Assessing pseudorandomness in numerical examples}\label{PRNGSect}

In this section, we describe a general method of empirical statistical assessment of any PRNG, and we apply it in sample examples to the PRNG $G_f$ associated to particular $2$-adic scaling maps $f:\ZZ_2\to\ZZ_2$.

A {\em pseudorandom number generator (PRNG)} is a triple $(G,M,N)$, where $M<N$ are positive integers and
\[
G:\{0,1\}^M\to \{0,1\}^N
\]
is any function.  The domain $\{0,1\}^M$ and codomain $\{0,1\}^N$ of $G$ are considered as probability spaces with the standard uniform measures $\PP_M$ and $\PP_N$, respectively.  The length $M$ bit strings in the domain $\{0,1\}^M$ of $G$ are viewed as {\em random seeds}.  Arbitrary length $N$ strings in the codomain $\{0,1\}^N$ of $G$ are {\em true random strings}, and length $N$ strings in the image $G(\{0,1\}^M)$ of $L$ are {\em pseudorandom strings}.  Intuitively, a PRNG is considered ``good'' if it is difficult for an observer to distinguish between pseudorandom strings and true random strings in $\{0,1\}^N$.

Let $\w\in\{0,1\}^r$ be an arbitrary binary word of length $r\geq1$, and let $X_\w$ be the number of times that the word $\w$ occurs in a random bit string in $\{0,1\}^N$.  As a random variable on $\{0,1\}^N$, $X_\w$ has mean $\mu_\w=(N+1-r)/2^r$ by linearity of expectation, since there are $N+1-r$ positions at which a length $r$ word can occur, and a probability of $1/2^r$ that $\w$ occurs at each such position.  Denote by $\sigma_\w$ the standard deviation of $X_\w$; there is a general formula for $\sigma_\w$ in terms of the set of self-compatible shifts of $\w$, but for our purposes it suffices to calculate this quantity in individual examples, as we do below.

For a parameter $q>0$ define 
\[
\mathscr{E}_{q,\w} = \{|X_\w-\mu_\w|\leq q\sigma_\w\}
\] 
to be the event in $\{0,1\}^N$ that the number of occurrences of the word $\w$ in a true random string is within $q$ standard deviations of its mean.    Consider also the $G$-inverse image of $\mathscr{E}_{q,\w}$, that is
\[
G^{-1}(\mathscr{E}_{q,\w})=\{\a\in\{0,1\}^M \mid G(\a)\in\mathscr{E}_{q,\w}\}.
\]
In other words, $G^{-1}(\mathscr{E}_{q,\w})$ is the event in the seed space $\{0,1\}^M$ that a random seed $\a$ generates a pseudorandom bit string $G(\a)$ satisfying $|X_\w(G(\a))-\mu_\w|\leq q\sigma_\w$.

We compare the probabilities $\PP_M(G^{-1}(\mathscr{E}_{q,\w}))$ and $\PP_N(\mathscr{E}_{q,\w})$.  If
\[
\PP_M(G^{-1}(\mathscr{E}_{q,\w}))\approx \PP_N(\mathscr{E}_{q,\w})\hskip1cm q=1,2,3,\dots,
\]
this indicates that the distribution of the frequency of the word $\w$ fails to meaningfully distinguish between pseudorandom and true random strings in $\{0,1\}^N$.  We may interpret this as one item of empirical evidence that the PRNG $G$ successfully produces low entropy, pseudorandom strings in $\{0,1\}^N$ which resemble high entropy, true random strings in $\{0,1\}^N$.  

We now apply this statistical assessment to the PRNG 
\[
G_f:\{0,1\}^M\to \{0,1\}^N
\]
associated to a $2$-adic scaling map $f:\ZZ_2\to\ZZ_2$.  We first recall the definition (\ref{ItineraryPRNGDef}) of $G_f$.  Given a bit string $a_0a_1a_2\dots a_{M-1}\in\{0,1\}^M$, assemble it into the $2$-adic expansion of a nonnegative integer in the range $0\leq \alpha\leq 2^M-1$ by 
\[
\alpha=a_0+a_12+a_22^2+\dots+a_{M-1}2^{M-1},
\]
and define
\begin{equation*}
G_f(a_0a_1a_2\dots a_{M-1}) = b_0b_1b_2\dots b_{N-1} \hskip1cm b_n=\overline{f^n(\alpha)}.
\end{equation*}

In Tables \ref{DistWordFreqTable1} and  \ref{DistWordFreqTable2}, we calculate $\PP_M(G_f^{-1}(\mathscr{E}_{q,\w}))$ and $ \PP_N(\mathscr{E}_{q,\w})$ in specific examples, taking $M=12$ and $N=240$ and thus representing PRNGs which increase the lengths of bit strings by a factor of 20.  Our data considers the cases $q=1,2,3$, all words $\w\in\{0,1,00,01,10,11\}$ of length at most $2$, for the particular $2$-adic scaling maps 
\begin{equation}
f(x) = \frac{1}{2}(x^2-x) \hskip1cm \text{and} \hskip1cm f(x) = \frac{3}{2}x^2-\frac{1}{2}x+1.
\end{equation} 
The data do seem to show a close resemblance between pseudorandom and true random word frequency distributions, supporting the viability of $G_f:\{0,1\}^M\to\{0,1\}^N$ as statistically sound PRNGs.

\begin{table}[ht]
\centering 
\caption{Word frequency distribution data for the PRNG $G_f:\{0,1\}^M\to \{0,1\}^N$ associated to the $2$-adic scaling map $f(x)=\frac{1}{2}(x^2-x)$, in the case $M=12$ and $N=240$ and all words $\w\in\{0,1,00,01,10,11\}$ of length $\leq 2$.}\label{DistWordFreqTable1}

\begin{tabular}{c|c|c|c}
$q$ & $\w$  & $\PP_M(G_f^{-1}(\mathscr{E}_{q,\w}))$ & $\PP_N(\mathscr{E}_{q,\w})$  \\
\hline
1 & 0,1 &  0.6630... & 0.6670... \\
2 & 0,1 &  0.9523... & 0.9548... \\
3 & 0,1 &  0.9958... & 0.9976... \\
1 & 00 &  0.6770... & 0.6740... \\
2 & 00 &  0.9550... & 0.9579... \\
3 & 00 &  0.9953... & 0.9971... \\
1 & 01  & 0.6967... & 0.6973... \\
2 & 01  & 0.9599... & 0.9609... \\
3 & 01  & 0.9948... & 0.9970... \\
1 & 10  & 0.6967... & 0.6973... \\
2 & 10  & 0.9589... & 0.9609... \\
3 & 10  & 0.9946... & 0.9970... \\
1 & 11 &  0.6667... & 0.6740... \\
2 & 11 &  0.9570... & 0.9579... \\
3 & 11 &  0.9956... & 0.9971... \\
\end{tabular}
\end{table}

\begin{table}[ht]
\centering 
\caption{Word frequency distribution data for the PRNG $G_f:\{0,1\}^M\to \{0,1\}^N$ associated to the $2$-adic scaling map $f(x)=\frac{3}{2}x^2-\frac{1}{2}x+1$, in the case $M=12$ and $N=120$ and all words $\w\in\{0,1,00,01,10,11\}$ of length $\leq 2$.}\label{DistWordFreqTable2}

\begin{tabular}{c|c|c|c}
$q$ & $\w$  & $\PP_M(G_f^{-1}(\mathscr{E}_{q,\w}))$ & $\PP_N(\mathscr{E}_{q,\w})$  \\
\hline
1 & 0,1 &  0.6555... & 0.6670... \\
2 & 0,1 &  0.9492... & 0.9548... \\
3 & 0,1 &  0.9978... & 0.9976... \\
1 & 00 &  0.6772... & 0.6740... \\
2 & 00 &  0.9533... & 0.9579... \\
3 & 00 &  0.9975... & 0.9971... \\
1 & 01  & 0.6828... & 0.6973... \\
2 & 01  & 0.9570... & 0.9609... \\
3 & 01  & 0.9963... & 0.9970... \\
1 & 10  & 0.6801... & 0.6973... \\
2 & 10  & 0.9599... & 0.9609... \\
3 & 10  & 0.9960... & 0.9970... \\
1 & 11 &  0.6589... & 0.6740... \\
2 & 11 &  0.9492... & 0.9579... \\
3 & 11 &  0.9975... & 0.9971... \\
\end{tabular}
\end{table}

The pseudorandom probabilities $\PP_M(G_f^{-1}(\Escr_{q,\w}))$ in Tables \ref{DistWordFreqTable1} and  \ref{DistWordFreqTable2} were calculated in MATLAB.  

The calculations of the true random probabilities $\PP_N(\mathscr{E}_{q,\w})$ rely on elementary combinatorial arguments, as follows.  For words $\w\in\{0,1\}$ of length $r=1$, the random variable $X_{\w}$ is binomial, with mean $\mu_\w=N/2$ and standard deviation $\sigma_\w=\sqrt{N}/2$, and 
\begin{equation*}
\begin{split}
\PP_N(\Escr_{q,\w})=\PP_N(|X_\w-\mu_\w|\leq q\sigma_\w)= \frac{1}{2^{N}}\sum_{A_{q,N}\leq n\leq B_{q,N}}\binom{N}{n}
\end{split}
\end{equation*}
where $A_{q,N} = \frac{N}{2}-q\frac{\sqrt{N}}{2}$ and $B_{q,N} = \frac{N}{2}+q\frac{\sqrt{N}}{2}$.

For words $\w\in\{00,01,10,11\}$ of length $r=2$, the calculation of $\PP_N(\mathscr{E}_{q,\w})$ is slightly more complicated, and we provide the full details here.  Recall the standard convention that the binomial coefficient $\binom{a}{b}$ takes the value zero when $b>a$.  

\begin{lem}\label{WordFreqLem}
Let $N\geq2$ and let $\w\in\{00,01,10,11\}$ be a binary word of length $2$.  Let $X_\w$ be the number of occurrences of the word $\w$ in a random string in $\{0,1\}^N$.  If $\w\in\{01,10\}$, then
\begin{equation}\label{FreqDist01Case}
\PP_N(X_{\w}=n)=\frac{1}{2^{N}}\binom{N+1}{2n+1}.
\end{equation}
If $\w\in\{00,11\}$, then
\begin{equation}\label{FreqDist11Case}
\PP_N(X_{\w}=n)=\frac{1}{2^{N}}\sum_{k=n+1}^{\lfloor(N+n+1)/2\rfloor}\binom{k-1}{n}\binom{N-k+1}{k-n}.
\end{equation}
\end{lem}

\begin{proof}
We first prove (\ref{FreqDist01Case}), and we treat the $\w=01$ case (the other case is identical by symmetry).  Clearly if $2n>N$ then 01 cannot occur $n$ times in a string of length $N$; thus both sides of (\ref{FreqDist01Case}) are zero.  

Assume now that $0\leq n\leq N/2$.  Let $S$ be any subset of $\{0,1,2,\dots,N\}$ with $|S|=2n+1$.  Define a bit string $\a=\a_S=a_0a_1a_2\dots a_{N-1}\in\{0,1\}^N$ as follows: 
\begin{itemize}
\item If $0\in S$ set $a_0=0$, and if $0\not\in S$ set $a_0=1$.
\item For $1\leq n\leq N-1$, recursively define $a_n\in\{0,1\}$ by $a_n=a_{n-1}$ if $n\not\in S$, and $a_n\neq a_{n-1}$ if $n\in S$.
\end{itemize}
If $0\in S$ then $a_0=0$ and there are $2n$ further indices in $S$ at which the string $\a$ oscillates from 0 to 1 or 1 to 0.  Exactly $n$ of these oscillations must be occurrences of 01, proving that $\a$ has exactly $n$ occurrences of 01.  If $0\not\in S$ then $a_0=1$ and there are $2n+1$ further indices in $S$ at which the string $\a$ oscillates from 0 to 1 or 1 to 0.  The first of these oscillations is an occurrence of 10, and of the remaining oscillations, exactly $n$ of these must be occurrences of 01, proving that $\a$ has exactly $n$ occurrences of 01 in this case as well.  

Conversely, any string $\a$ in $\{0,1\}^N$ with exactly $n$ occurrences of 01 uniquely determines a subset $S$ of $\{0,1,2,\dots,N\}$ with $|S|=2n+1$ and $\a=\a_S$, by declaring that $0\in S$ if and only if $a_0=0$, and for $n\geq1$ declaring that $n\in S$ if and only if $a_n\neq a_{n-1}$.  

The preceding discussion has shown that $S\mapsto \a_S$ is a bijection from the set of all subsets of $\{0,1,2,\dots,N\}$ with $|S|=2n+1$ and the set of all strings in $\{0,1\}^N$ with exactly $n$ occurrences of 01.  This bijection proves the identity (\ref{FreqDist01Case}).

We turn now to the proof of (\ref{FreqDist11Case}) for $\w\in\{00,11\}$, and we treat the $\w=11$ case (the other case is identical by symmetry).  We count the number of strings $\a\in\{0,1\}^N$ having exactly $n$ occurrences of 11.  Let $\a$ be a such a string and let $k$ be the total number of 1's occurring in the string.  Thus $k\geq n+1$, since for any string containing $k$ 1's, the most number of occurrences of $11$ is $k-1$, which is possible only when all $k$ 1's are consecutive.  

The number of runs of 1s in the string $\a$ is $r=k-n$, because the number of occurrences of 11 in any run of length $\ell$ is $\ell-1$.  For fixed $n$ and $k$ as above, the string $\a$ is now uniquely determined by two choices: the way that the $k$ 1s in the string are divided into $r$ runs, and the way that these $r$ runs are distributed among the $N$ total positions of the string $\a$.

The number of ways to divide a collection of $k$ 1s, say $1_{1}1_{2}1_{3}\dots1_{k}$, into $r$ runs is 
\[
\binom{k-1}{r-1}=\binom{k-1}{k-1-n}=\binom{k-1}{n}
\]
since each such configuration is equivalent to the choice of $r-1$ indices in $1,2,3,\dots,k-1$ to be the final indices in the first $r-1$ of these runs.  The number of ways to distribute these $r$ runs among the $N$ positions of the string $\a$ is 
\[
\binom{N-k+1}{r}=\binom{N-k+1}{k-n}
\]
because the complement of the 1s are the $N-k$ 0s, say $0_10_20_2\dots0_{N-k}$, and each configuration of the $r$ runs of 1s within these 0's is equivalent to the choice of $r$ out of of $N-k+1$ available slots in between consecutive $0_i's$, including the slot preceding $0_1$ and the slot following $0_{N-k}$.  Note that when $k>(N+n+1)/2$, it follows that $r=k-n>N-k+1$, and no such configuration is possible; hence no strings $\a$ of length $N$ exist with $k>(N+n+1)/2$ total 1's and $n$ occurrences of $11$.

We conclude from these two counts that the total number of strings $\a\in\{0,1\}^N$ having $k$ 1's and exactly $n$ occurrences of 11 is 
\[
\binom{k-1}{n}\binom{N-k+1}{k-n},
\]
and summing over all possible values $n+1\leq k\leq \lfloor(N+n+1)/2\rfloor$, the calculation (\ref{FreqDist11Case}) follows.
\end{proof}

\begin{prop}
Let $N\geq2$ and let $\w\in\{00,01,10,11\}$ be a binary word of length $2$.  Let $X_\w$ be the number of occurrences of the word $\w$ in a random string in $\{0,1\}^N$.  

\noindent
{\bf (a)} As a random variable, $X_\w$ has mean  $\mu_\w=(N-1)/4$ and standard deviation 
\begin{equation*}
\sigma_{\w}= 
\begin{cases}
\frac{\sqrt{5N-7}}{4} & \text{ if } \w\in\{00,11\} \\
\frac{\sqrt{N+1}}{4} & \text{ if } \w\in\{01,10\}.
\end{cases}
\end{equation*}

\noindent
{\bf (b)} If $w\in\{01,10\}$, then 
\begin{equation*}
\begin{split}
\PP_N(\Escr_{q,\w})= \frac{1}{2^{N}}\sum_{A_{q,N}\leq n\leq B_{q,N}}\binom{N+1}{2n+1}
\end{split}
\end{equation*}
where $A_{q,N} = \frac{N-1}{4}-q\frac{\sqrt{N+1}}{4}$ and $B_{q,N} = \frac{N-1}{4}+q\frac{\sqrt{N+1}}{4}$.  

\noindent
{\bf (c)} If $w\in\{00,11\}$, then
\begin{equation*}
\begin{split}
\PP_N(\Escr_{q,\w}) & = \frac{1}{2^{N}}\sum_{A_{q,N}\leq n\leq B_{q,N}}\sum_{k=n+1}^{\lfloor(N+n+1)/2\rfloor}\binom{k-1}{n}\binom{N-k+1}{k-n}
\end{split}
\end{equation*}
where $A_{q,N} = \frac{N-1}{4}-q\frac{\sqrt{5N-7}}{4}$ and $B_{q,N} = \frac{N-1}{4}+q\frac{\sqrt{5N-7}}{4}$.
\end{prop}
\begin{proof}
The proof of {\bf (a)} is a straightforward calculation obtained by writing $X_\w=Y_0+Y_1+Y_2+\dots+Y_{N-2}$, where $Y_i=1$ if $\w$ occurs in the $i$-th position of a random string in $\{0,1\}^N$, and $Y_i=0$ otherwise.  Parts {\bf (b)} and {\bf (c)} follow from Lemma \ref{WordFreqLem} and the definition $\Escr_{q,\w}=\{|X_\w-\mu_\w|\leq q\sigma_\w\}$.
\end{proof}

\def\cprime{$'$}



\begin{thebibliography}{10}

\bibitem{Andrecut}
{\sc M.~Andrecut}, {\em On combination of stochastic and deterministic
  processes}, Int. J. of Modern Phys., 12 (1998), pp.~921--930.

\bibitem{AntunesHillPRNGML}
{\sc B.~A. Antunes and D.~Hill}, {\em Random numbers for machine learning},
  Journal of Data Science and Intelligent Systems, 2024.

\bibitem{MR3890051}
{\sc R.~L. Benedetto}, {\em Dynamics in one non-archimedean variable}, vol.~198
  of Graduate Studies in Mathematics, American Mathematical Society,
  Providence, RI, 2019.

\bibitem{dieharder}
{\sc R.~G. Brown, D.~Eddelbuettel, and D.~Bauer}, {\em Dieharder: a random
  number test suite},
  \href{https://webhome.phy.duke.edu/~rgb/General/dieharder.php}{https://webhome.phy.duke.edu/$\sim$rgb/General/dieharder.php}.

\bibitem{KatzLindell}
{\sc J.~Katz and Y.~Lindell}, {\em Introduction to {M}odern {C}ryptography:
  {P}rinciples and {P}rotocols}, CRC Press, 2007.

\bibitem{MR2439398}
{\sc J.~Kingsbery, A.~Levin, A.~Preygel, and C.~E. Silva}, {\em On
  measure-preserving {$C\sp 1$} transformations of compact-open subsets of
  non-{A}rchimedean local fields}, Trans. Amer. Math. Soc., 361 (2009),
  pp.~61--85.

\bibitem{MR2773138}
\leavevmode\vrule height 2pt depth -1.6pt width 23pt, {\em Dynamics of the
  {$p$}-adic shift and applications}, Discrete Contin. Dyn. Syst., 30 (2011),
  pp.~209--218.

\bibitem{MR1556918}
{\sc J.~F. Koksma}, {\em Ein mengentheoretischer {S}atz \"uber die
  {G}leichverteilung modulo {E}ins}, Compositio Math., 2 (1935), pp.~250--258.

\bibitem{MachicaoBruno}
{\sc J.~Machicao and O.~M. Bruno}, {\em Improving the pseudo-randomness
  properties of chaotic maps using deep-zoom}, Chaos., 27 (2017).

\bibitem{MR836749}
{\sc H.~Niederreiter}, {\em Multidimensional numerical integration using
  pseudorandom numbers}, no.~27, 1986, pp.~17--38.
\newblock Stochastic programming 84. I.

\bibitem{PareekPatidarSud}
{\sc N.~K. Pareek, V.~Patidar, and K.~K. Sud}, {\em A random bit generator
  using chaotic maps}, Int. J. of Network Security, 10 (2010), pp.~32--38.

\bibitem{MR3958063}
{\sc C.~Petsche and E.~Stacy}, {\em A dynamical construction of small totally
  {$p$}-adic algebraic numbers}, J. Number Theory, 202 (2019), pp.~27--36.

\bibitem{PhatakRao}
{\sc S.~C. Phatak and S.~S. Rao}, {\em Logistic map: A possible random-number
  generator}, Phys. Rev. E, 51 (1995), p.~3670.

\bibitem{RosulekTJOC}
{\sc M.~Rosulek}, {\em The Joy of Cryptography}, MIT Press, Cambridge, MA,
  2025.

\bibitem{NIST}
{\sc A.~Rukhin, J.~Soto, J.~Nechvatal, M.~Smid, E.~Barker, S.~Leigh,
  M.~Levenson, M.~Vangel, D.~Banks, A.~Heckert, J.~Dray, and S.~Vo}, {\em A
  statistical test suite for random and pseudorandom number generators for
  cryptographic applications}, National Institute of Standards and Technology,
  (2010).

\bibitem{MR2316407}
{\sc J.~H. Silverman}, {\em The {A}rithmetic of {D}ynamical {S}ystems},
  vol.~241 of Graduate Texts in Mathematics, Springer, New York, 2007.

\bibitem{UlamVonNeumann}
{\sc S.~M. Ulam and J.~von Neumann}, {\em On combination of stochastic and
  deterministic processes}, Bull. Amer. Math. Soc., 53 (1947).

\bibitem{woodcocksmart:chaosprng}
{\sc C.~F. Woodcock and N.~P. Smart}, {\em $p$-adic chaos and random number
  generation}, Experiment. Math., 7 (1998), pp.~333--342.

\end{thebibliography}

\end{document}